\documentclass[12pt]{article}
\usepackage{url}
\renewcommand{\baselinestretch}{1.150}

\usepackage{amsmath}
\usepackage{amsthm}
\usepackage{amssymb, amsfonts}
\usepackage{mathrsfs}
\newcommand{\ds}{\displaystyle}
\newcommand{\p}{\partial}

\newtheorem{rmk}{Remark}

\newcommand{\RO}{recursion operator}

\DeclareMathAlphabet{\bi}{OML}{cmm}{b}{it}

\newcommand{\be}{\begin{equation}}
\newcommand{\ee}{\end{equation}}
\newcommand{\ba}{\begin{array}}
\newcommand{\ea}{\end{array}}

\date{October 15, 2017}

\title{\protect\vspace{-2cm}\sc\Large Recursion Operators\\ for Multidimensional Integrable PDEs\protect\vspace{-0.3cm}}
\author{{\sc A. Sergyeyev}\\[5mm]
Mathematical Institute, Silesian University in Opava,\\ Na
Rybn\'\i{}\v{c}ku 1, 746 01 Opava,~Czech Republic\\
E-mail: \tt{Artur.Sergyeyev@math.slu.cz}}

\begin{document}
\maketitle

\begin{abstract}\protect\vspace*{-12mm}
We present a novel construction of recursion operators for scalar second-order integrable multidimensional PDEs with isospectral Lax pairs written in terms of first-order scalar differential operators. Our approach 
is quite straightforward and can be readily applied using modern computer algebra software. It is illustrated by examples, two of which are new. \looseness=-1
\vspace{-5mm}
\end{abstract}


\section{Introduction}
The search for integrable partial differential equations in multidimensions, i.e., in three or more independent variables, is among the most important problems in the modern theory of integrable systems, cf.\ e.g.\ \cite{ac,bls,d-book,kvv,as,aslmp} and references therein. It is well known that such PDEs are never isolated -- they belong to infinite countable families of pairwise compatible integrable systems; such families are known as integrable hierarchies. Existence of such an hierarchy plays an important role in the study of PDE in question, including the construction of its exact solutions, cf.\ e.g.\ \cite{ac,d-book,kvv,o} and references therein. The construction of an integrable hierarchy associated with a given PDE is usually most easily performed using a recursion operator which, roughly speaking, maps any given symmetry of the PDE under study into a (new) symmetry, see e.g.\ \cite{o} and references therein.

In what follows we restrict ourselves to considering a second-order partial differential equation in $d$ independent variables $x^1,\dots,x^d$ for a single unknown function $u$,
\begin{equation}\label{e}
F(\vec x,u,u_{x^1},\dots,u_{x^d}, u_{x^1 x^1},\dots, u_{x^d x^d})=0.
\end{equation}
Here $\vec x=(x^1,\dots,x^d)$, and, as usual, the subscripts refer to partial derivatives. All functions here and below are assumed sufficiently smooth for all computations to make sense; this can be formalized using the language of differential algebra, cf.\ e.g.\ the discussion in \cite{aslmp}.

This class contains plenty of integrable equations, e.g.\ the dispersionless Kadomtsev--Petviashvili also known as Khokhlov--Zabolotskaya equation \cite{kz}, the Pavlov equation \cite{pav}, various heavenly equations \cite{df}, the Boyer--Finley equation \cite{bf}, the ABC equation \cite{ksm, z} etc., and many of those equations are of relevance for applications.

To the best of our knowledge, for nonlinear equations of the form (\ref{e}) in three or more independent variables all \RO{}s known to date are in fact B\"acklund auto-transformations for the linearized versions of these equations rather than genuine operators mapping one symmetry to another, see e.g.\ \cite{kvv,ms,p,as} for details and references. 

To date, several methods for construction of such recursion operators  were discovered,
see \cite{m-s,ms,m1,m2,m3,as} and references therein for details. However, most of these methods either involve rather heavy computations (those of \cite{ms}, \cite{m1, m2, m3}, and, to a somewhat lesser extent, \cite{as}) or have a somewhat resricted range of applicability (that of \cite{m-s}).

Below we introduce an approach to the construction of recursion operators in question which, while inspired by the ideas from \cite{as}, is significantly more straightforward and not too computationally demanding yet apparently not really narrower than that of \cite{as} in the range of applicability as far as equations of the form (\ref{e}) are concerned.

The rest of the paper is organized as follows. In Section~\ref{pre} we quickly recall the basic definitions and facts needed for the sequel, and in Section~\ref{c} we present our approach to construction of the recursion operators for equations of the form (\ref{e}). Section~\ref{exa} contains examples, and Section~\ref{d} gives conclusions and discussion.
\section{Preliminaries}\label{pre}

We shall need the operator of linearization of the left-hand side $F$ of (\ref{e}), cf.\ e.g.\ \cite{kvv, kv},
\[
\ell_F=\displaystyle\frac{\p F}{\p u}+
\sum_{i=1}^{d}\frac{\p F}{\p u_{x^i}} D_{x^i}+\sum_{i=1}^{d}\sum_{j=i}^d
\frac{\p F}{\p u_{x^i x^j}} D_{x^i} D_{x^j},
\]
where $D_{x^i}$ are total derivatives, see e.g.\ \cite{kvv,kv,o,as} and references therein for details on those.

For the rest of a paper we make a {\em blanket assumption} that 
(\ref{e}) 
admits a 
Lax pair of the form
\begin{equation}\label{lax-x}
\mathscr{X}_i\psi=0,\quad i=1,2,
\end{equation}
where $\mathscr{X}_i$ are linear
in the spectral parameter $\lambda$ first-order 
differential operators of the form
\[
\mathscr{X}_i=\lambda\mathscr{X}_i^1-\mathscr{X}_i^0,
\]
where
\[
\mathscr{X}_i^s=X_i^{s,0}+\sum\limits_{j=1}^d X_i^{s,j} D_{x^j},\ i=1,2, \ s=0,1,
\]
and $X_i^{s,k}$ for all $i,s,k$ 
depend on $\vec x$, $u$ and first- and possibly second-order derivatives of $u$; see e.g.\ \cite{m-s,zs} on Lax pairs of this type.

\begin{rmk}\em A Lax pair of the form (\ref{lax-x}), especially with $X_i^{s,0}\equiv 0$, typically exists only if the PDE (\ref{e}) under study is dispersionless, i.e., it can be written as a first-order homogeneous quasilinear system, cf.\ e.g.\ \cite{al, bls, ck, dfkn, as} and references therein for details.

This is the case e.g.\ for any quasilinear PDE of the form
\begin{equation}\label{dsl}
\sum_{i=1}^d \sum_{j=i}^d f^{ij}(\vec x, u_{x^1},\dots,u_{x^d})u_{x^i x^j}=0,
\end{equation}
which can (see e.g.\ \cite{as} for details) be turned into a quasilinear homogeneous first-order system for
the vector function
$\boldsymbol{v}=(u_{x^1}, \dots, u_{x^d})$
upon being supplemented by the natural compatibility conditions
\[
(v_i)_{x^j}=(v_j)_{x^i},\quad i,j=1,\dots,d,\quad i<j.
\]

Of course, 
integrable dispersionless second-order PDEs are rather exceptional, 
and some of them admit Lax pairs of a significantly more general form than (\ref{lax-x}), for instance,
involving higher powers of $\lambda$ and/or derivatives w.r.t.\ $\lambda$, cf.\ e.g.\ the discussion in \cite{ck, m-s, as, aslmp} and references therein.

Note, however, that in certain cases the Lax pairs involving second power of $\lambda$ (but not the derivatives w.r.t.\ $\lambda$) can still be rewritten in the form which is linear in $\lambda$, cf.\ e.g.\ Example 3 below and the treatment of the Pavlov equation in \cite{as}. 
\end{rmk}


\section{The construction of recursion operators}\label{c}

Now return to (\ref{e}) and (\ref{lax-x}) and suppose for a moment that $\psi$ defined by (\ref{lax-x}) also satisfies
\[
\ell_{{F}}(\psi)=0
\]
by virtue of (\ref{e}), (\ref{lax-x}) and differential consequences thereof, i.e., $\psi$ is
a nonlocal symmetry\footnote{For the sake of brevity in what follows we shall refer to solutions of $\ell_F(U)=0$ just as to nonlocal symmetries rather than shadows (the latter terminology is employed e.g.\ in \cite{kvv,kv} and references therein).} for (\ref{e}).

Then substituting a formal Laurent expansion $\psi=\ds\sum\limits_{k=0}^\infty \psi_k \lambda^{-k}$ into (\ref{lax-x}) 
gives rise (cf.\ e.g.\ \cite{p, as} and references therein)
to the recursion relations
\[
\mathscr{X}_i^1(\psi_{k+1})=\mathscr{X}_i^0 (\psi_{k}),\ i=1,2,\quad k=0,1,2,\dots
\]
which define 
an infinite hierarchy of nonlocal symmetries $\psi_k$,  $k=0,1,2,\dots$, for (\ref{e}).

This suggests, see \cite{p,as}, that more general relations
\begin{equation}\label{ro0}
\mathscr{X}_i^1(\tilde U)=\mathscr{X}_i^0 (U),\ i=1,2
\end{equation}
could define a recursion operator for (\ref{e}) which maps an arbitrary (nonlocal) symmetry $U$ of (\ref{e}) to a new (again in general nonlocal) symmetry $\tilde U$ of (\ref{e}); checking whether (\ref{ro0}) indeed defines a recursion operator for (\ref{e}) is a straightforward matter.

The situation when $\psi$ from the original Lax pair is a nonlocal symmetry and one can construct the recursion operator in the above fashion occurs e.g.\ for the general heavenly equation and for the six-dimensional heavenly equation, whose recursion operators were found respectively in \cite{as} and in \cite{ms}.\looseness=-1


However, in general the relations (\ref{ro0}) constructed from the original Lax pair (\ref{lax-x}) do not define the recursion operator.
%
In this case we propose to twist these relations, namely, to replace (\ref{ro0}) by
\begin{equation}\label{rot}
\mathscr{X}_i^1(\tilde U)+f_i^1\tilde U=f_i^0 U+\mathscr{X}_i^0 (U),\quad i=1,2,
\end{equation} 
where  $f_i^s$, just like $X_i^{s,j}$ above, 
depend on $\vec x$, $u$ and first- and possibly second-order derivatives of $u$,
and require that relations (\ref{rot}) define a recursion operator for (\ref{e}).

In other words, one should require that

\medskip

\begin{minipage}{0.9\textwidth}

\begin{enumerate}\itemsep=-1mm
\item[a)] system (\ref{rot}), seen as a system for $\tilde U$, is compatible
if $U$ is a symmetry for (\ref{e});

\item[b)] $\tilde U$ defined by (\ref{rot}) is a (in general nonlocal) symmetry for (\ref{e}) if so is $U$.

\end{enumerate}
\end{minipage}

\medskip


Imposing the requirements in question gives rise to a (in general overdetermined) system of equations for $f_i^j$, and any solution of this system gives a \RO{} for (\ref{e}). 

Writing down and solving the said system of equations for $f_i^j$ is quite straightforward and can be readily done using the existing computer algebra software,
for instance \cite{m} or \cite{v,v1}.

Thus, we propose the following procedure for searching for a recursion operator for (\ref{e}):

\medskip

\begin{minipage}{0.9\textwidth}

\begin{enumerate}\itemsep=-1mm
\item Write down the relations (\ref{rot}) and require that the conditions a) and b) are satisfied.

\item To solve the system of equations for $f_i^j$ resulting from the previous step.

\end{enumerate}
\end{minipage}

\medskip

Any solution of the said system for $f_i^j$ gives rise to a recursion operator for (\ref{e}) defined via (\ref{rot}).

In principle such solutions are by no means obliged to exist but we yet have to find a case where such nonexistence really occurs (bearing in mind, of course, our blanket assumption, i.e., that (\ref{e}) has a nontrivial Lax pair of the form (\ref{lax-x})).\looseness=-1 

\section{Examples}\label{exa}
\noindent{\bf Example 1.}  Consider the six-dimensional equation
\be\label{equ}
u_s u_{zt}-u_z u_{st}- u_s u_{xy}+u_y u_{sx}-u_y u_{rz}+u_z u_{ry}=0,
\ee
derived in \cite{as} from a certain first-order system from \cite{fk}.

This equation has \cite{as} a Lax pair of the form (\ref{lax-x}) with 
\[
\mathcal{X}_1=\displaystyle D_{z}-\frac{u_z}{u_s}D_{s}-\lambda D_{x}+\lambda \frac{u_z}{u_s} D_{r},\quad
\mathcal{X}_2=\displaystyle D_{y}-\frac{u_y}{u_s}D_{s}-\lambda D_{t}+\lambda \frac{u_y}{u_s}D_{r}.
\]

The \RO{} 
for (\ref{equ}) found in \cite{as}, which enables us to construct a (nonlocal) symmetry $\tilde U$ of (\ref{equ}) from a (nonlocal) symmetry $U$ of (\ref{equ}), can be written as 
\be\label{ro1}
\ba{l}
\tilde{U}_y-\ds \frac{u_y}{u_s} \tilde{U}_s=-\frac{u_y}{u_s} U_r+ U_t-\frac{(u_{st}-u_{ry})}{u_s} U, \\[7mm]
\tilde{U}_z-\ds \frac{u_z}{u_s} \tilde{U}_s=-\frac{u_z}{u_s} U_r+ U_x+\frac{(u_{rz}-u_{sx})}{u_s} U,
\ea
\ee
and so obviously is of the form (\ref{rot}) with
\[
f_1^0=-\ds\frac{(u_{st}-u_{ry})}{u_s},\quad f_2^0=\ds\frac{(u_{rz}-u_{sx})}{u_s},\quad f_1^1=0,\quad f_2^1=0.
\]

\noindent{\bf Example 2.}  Consider the equation \cite{dfkn}
\begin{equation}\label{ex2}
D_z(u_y/u_x)-D_x(u_t/u_x)=0
\end{equation}
which has a Lax pair \cite{dfkn}
\[
\psi_t=\lambda\psi_z+(u_t/u_x)\psi_x,\quad \psi_y=(\lambda+u_y/u_x)\psi_x,
\]
so we have $\mathcal{X}_i^{s,0}=0$ for all $i$ and $s$, $\mathcal{X}_1^0=D_t-(u_t/u_x)D_x$, $\mathcal{X}_2^0=D_y-(u_y/u_x)D_x$, $\mathcal{X}_1^1=D_z$, $\mathcal{X}_2^1=D_x$.

A straightforward computation along the lines of previous section shows that a \RO{} for (\ref{ex2}) of the form (\ref{rot}) indeed can be constructed using our method presented above.

Namely, it is easily seen that for
\[
f_1^1=\ds-\frac{u_{xz}}{u_x},
\quad f_2^1=\ds-\frac{u_{xx}}{u_x},
\quad f_1^0=0,\quad f_2^0=0,
\]
and $\mathcal{X}_i^j$ given above
the relations
\be\label{ro-dfkn}
\ba{l}
\tilde{U}_z\ds-\frac{u_{xz}}{u_x} \tilde{U}=U_t-\ds \frac{u_t}{u_x} {U}_x, \\[7mm]
\tilde{U}_x- \ds\frac{u_{xx}}{u_x} \tilde{U} =U_y-\ds \frac{u_y}{u_x} {U}_x.
\ea
\ee
define a \RO{} for (\ref{ex2})  which, to the best of our knowledge, has not yet appeared in the literature. 

\noindent{\bf Example 3.} Consider an equation of the form \cite{dfkn}
\begin{equation}\label{ex3}
D_t(m)+\alpha m D_x(n)-D_z(n)-\alpha n D_x(m)=0
\end{equation}
where $m=(u_y-u_z)/u_x$ and $n=(u_z-u_t)/u_x$, and $\alpha$ is an arbitrary constant.

It has \cite{dfkn} a Lax pair of the form
\[
\psi_y=c (m+\lambda n)\psi_x+\lambda^2\psi_t,\quad \psi_z=c n\psi_x+\lambda\psi_t
\]
where $c=1+\alpha-\lambda\alpha$.

The above Lax pair is quadratic in $\lambda$ but it can be rewritten in the form which is already linear in $\lambda$, i.e., which belongs to the class (\ref{lax-x}):
\[
\psi_y=\lambda \psi_z+c m\psi_x,\quad \psi_z=c n\psi_x+\lambda\psi_t.
\]

Applying the procedure from the preceding section to the above Lax pair 
immediately gives rise to an apparently hitherto unknown recursion operator for (\ref{ex3}) of the form (\ref{rot}):
\be\label{ro-dfkn-2}
\begin{array}{l}
\tilde U_z-\displaystyle\alpha\frac{(u_y-u_z)}{u_x}\tilde U_x+f_1^1 \tilde U=-(1+\alpha)\frac{(u_y-u_z)}{u_x}U_x+U_y+f_1^0 U,\\[5mm]
\tilde U_t+\displaystyle\alpha\frac{(u_t-u_z)}{u_x}\tilde U_x+f_2^1 \tilde U=\displaystyle(1+\alpha)\frac{(u_t-u_z)}{u_x}U_x+U_z+f_2^0 U,
\end{array}
\ee
where
\[
f_1^0=f_1^1=\frac{(\alpha(u_{xy}-u_{xz})-u_{xz})}{u_x}, \quad f_2^0=f_2^1=\frac{(\alpha(u_{xz}-u_{xt})-u_{xt})}{u_x}.
\]
\section{Conclusions and discussion}\label{d}
We have presented above a novel construction for recursion operators of multidimensional integrable second-order scalar PDEs possessing isospectral Lax pairs written in terms of first-order differential operators and linear in the spectral parameter $\lambda$, and illustrated this construction by examples, two of which are, to the best of our knowledge, new.

It appears that using our construction it is possible, in spite of its simplicity, to reproduce all known to date recursion operators for multidimensional equations of the form (\ref{e}), including e.g.\ those for the Pavlov equation \cite{pav}, six-dimensional heavenly equation, see \cite{ms} for its recursion operator, all integrable four-dimensional symplectic Monge--Amp\`ere equations listed in \cite{df},
the ABC equation with $A+B+C=0$, see \cite{ms} for its \RO{} and also \cite{ksm}, the Mart\'\i{}nez Alonso--Shabat equation \cite{m3}, etc. In fact, it would be very interesting to find a counterexample, i.e., a recursion operator for an equation of the form (\ref{e}) that would {\em not} be reproducible within our approach, in other words, that cannot be written in the form (\ref{rot}) for suitable $\mathcal{X}_i^{j}$ and $f_i^j$, but so far we failed to do so.\looseness=-1

In closing note that our method for construction of \RO{}s has no inherent restriction on the order of the equation under study: in principle, it is applicable not just to second-order PDEs (\ref{e}) but also to more general equations of the form
\begin{equation}\label{k}
F(\vec x,u,\overset{(k)}{u})=0,
\end{equation}
where $\overset{(k)}{u}$ denotes the derivatives of $u$ up to order $k$, and $k\geqslant 2$.

The key requirement for applicability of our method in this extended setup is that (\ref{k}) admits a Lax pair of the form (\ref{lax-x}) where now $X_i^{s,j}$ can depend on $\vec x$, $u$, and $\overset{(k)}{u}$, and the same goes for $f_i^j$.
\subsection*{Acknowledgments}
The author is greatly pleased to thank R. Vitolo for the inspiring discussions.\looseness=-1

This research was supported in part by the Ministry of Education, Youth and Sports of the Czech Republic (M\v{S}MT \v{C}R) under RVO funding for I\v{C}47813059, and by the Grant Agency of the Czech Republic (GA \v{C}R) under grant P201/12/G028.


\end{document}